\documentclass{amsart}%
\usepackage{amssymb}
\usepackage{amsmath}
\usepackage{latexsym}
\usepackage{hyperref}
\usepackage{amsfonts}
\usepackage{graphicx}%
\newcommand{\be}{\begin{equation}}
\newcommand{\ee}{\end{equation}}
\newcommand{\bea}{\begin{eqnarray}}
\newcommand{\eea}{\end{eqnarray}}

\begin{document}
\title{An Elementary Remark on a Lower bound for the $\mu$-Invariant of Singularity
Models for Ricci Flow}
\author{Bennett Chow$^{1}$}
\maketitle

This is a purely expository note.\footnotetext[1]{Address: Department of
Mathematics, University of\ California San\ Diego. La Jolla, CA 92093} Assume
the $C^{\infty}$ category; convergence of metrics or solutions is in the
pointed (Hamilton--)Cheeger--Gromov sense. Given a complete Riemannian
manifold $\left(  \mathcal{N}^{n},g\right)  $, let $\mathcal{W}\left(
g,f,\tau\right)  $ denote Perelman's entropy functional, where $f:\mathcal{N}%
\rightarrow\mathbb{R}$ and $\tau>0$; let $\mu\left(  g,\tau\right)  $ be its
infimum with the usual constraint and with $e^{-f/2}$ having compact support.
Recall that under metric convergence $\left(  \mathcal{M}_{k}^{n}%
,g_{k}\right)  \rightarrow\left(  \mathcal{M}_{\infty}^{n},g_{\infty}\right)
$, $\mu\left(  g_{\infty},\tau\right)  \geq\limsup_{k\rightarrow\infty}%
\mu\left(  g_{k},\tau\right)  $, for $\tau>0$.

Let $\left(  \mathcal{M}^{n},g\left(  t\right)  \right)  $, $0\leq t<T<\infty
$, be a singular solution to the Ricci flow on a closed manifold. Perelman's
monotonicity formula \cite{Perelman1} implies $\mu\left(  g\left(  t\right)
,\tau\right)  \geq\mu\left(  g\left(  0\right)  ,\tau+t\right)  $. Suppose for
$x_{i}\in\mathcal{M}$, $t_{i}\nearrow T$, and $K_{i}\rightarrow\infty$ that
$\left(  \mathcal{M},g_{i}\left(  t\right)  ,\left(  x_{i},0\right)  \right)
$, $g_{i}\!\left(  t\right)  \!=\!K_{i}g\!\left(  t_{i}+K_{i}^{-1}t\right)  $,
converges to a complete noncompact nonflat ancient solution $\left(
\mathcal{M}_{\infty}^{n},g_{\infty}\left(  t\right)  ,\left(  x_{\infty
},0\right)  \right)  $, $t\leq0$ (a singularity model). At any time, it is
conceivable that the curvature of $g_{\infty}\left(  t\right)  $ is unbounded.
By the limit, scaling ($\mu\left(  cg,c\tau\right)  =\mu\left(  g,\tau\right)
$), monotonicity, and continuity properties of $\mu$, for $t\leq0$, $\tau>0$,
and $\bar{t}<T$,%
\begin{align*}
\mu\left(  g_{\infty}\left(  t\right)  ,\tau\right)   &  \geq\limsup
_{i\rightarrow\infty}\mu\left(  g\left(  t_{i}+K_{i}^{-1}t\right)  ,K_{i}%
^{-1}\tau\right)  \\
&  \geq\limsup_{i\rightarrow\infty}\mu\left(  g\left(  \bar{t}\,\right)
,t_{i}-\bar{t}+K_{i}^{-1}\left(  \tau+t\right)  \right)  \\
&  =\mu\left(  g\left(  \bar{t}\,\right)  ,T-\bar{t}\,\right)  .
\end{align*}
Hence, any singularity model has the property that, for $t\leq0$ and $\tau>0$,
we have%
\[
\mu\left(  g_{\infty}\left(  t\right)  ,\tau\right)  \geq\lim_{\bar
{t}\rightarrow T}\mu\left(  g\left(  \bar{t}\,\right)  ,T-\bar{t}\,\right)
\doteqdot\mu_{T}\geq\mu\left(  g\left(  0\right)  ,T\right)  .
\]
Recall, by Perelman's no local collapsing result, there exists $\kappa>0$ such
that for any $\left(  x,t\right)  $ and $r\leq1$, if $R_{g_{\infty}\left(
t\right)  }\leq r^{-2}$ in $B_{g_{\infty}\left(  t\right)  }\left(
x,r\right)  $, then $\operatorname{Vol}_{g_{\infty}\left(  t\right)
}B_{g_{\infty}\left(  t\right)  }\left(  x,r\right)  \geq\kappa r^{n}$. In
particular, if $R_{g_{\infty}\left(  t\right)  }\leq1$ (such as for a
normalized nonflat steady gradient Ricci soliton singularity model), then
$\operatorname{Vol}_{g_{\infty}\left(  t\right)  }B_{g_{\infty}\left(
t\right)  }\left(  x,r\right)  \geq\kappa r^{n}$ for all $x\in\mathcal{M}%
_{\infty}$ and $r\leq1$. Moreover, on a unit ball $B\left(  x,1\right)  $ in
such a steady singularity model, Bakry--Emery volume comparison is comparable
to Riemannian volume comparison since we may add a constant to $f$ to make
$f\left(  x\right)  =0$, so that $\left\vert f\left(  y\right)  \right\vert
\leq1$ for $y\in B\left(  x,1\right)  $ by $\left\vert \nabla f\right\vert
\leq1$. As observed by O.\ Munteanu and J.\ Wang \cite{MunteanuWang}, this
implies that for $B\left(  x,1\right)  $ we have uniform bounds (independent
of $x$) for the constants of Cheeger, the Neumann Poincar\'{e} inequality, a
Sobolev inequality, and an isoperimetric inequality.

\end{document}